# Complex Datasets and Inverse Problems

## Tomography, Networks and Beyond

Regina Liu, William Strawderman and Cun-Hui Zhang, Editors

# Contents

# Preface

This book is a collection of papers dedicated to the memory of Yehuda Vardi. Yehuda was the chair of the Department of Statistics of Rutgers University when he passed away unexpectedly on January 13, 2005. On October 21–22, 2005, some 150 leading scholars from many different fields, including statistics, telecommunications, biomedical engineering, bioinformatics, biostatistics and epidemiology, gathered at Rutgers in a conference in his honor. This conference was on "Complex Datasets and Inverse Problems: Tomography, Networks, and Beyond," and was organized by the editors. The present collection includes research work presented at the conference, as well as contributions from Yehuda's colleagues.

The theme of the conference was networks and other important and emerging areas of research involving incomplete data and statistical inverse problems. Networks are abundant around us: communication, computer, traffic, social and energy are just a few examples. As enormous amounts of network data are collected in this information age, the field has attracted a great amount of attention from researchers in statistics and computer engineering as well as telecommunication providers and various government agencies. However, few statistical tools have been developed for analyzing network data as they are typically governed by time-varying and mutually dependent communication protocols sitting on complicated graph-structured network topologies. Many prototypical applications in these and other important technologies can be viewed as statistical inverse problems with complex, massive, high-dimensional and possibly biased/incomplete data. This unifying theme of inverse problems is particularly appropriate for a conference and volume dedicated to the memory of Yehuda. Indeed he made influential contributions to these fields, especially in medical tomography, biased data, statistical inverse problems, and network tomography.

The conference was supported by the NSF Grant DMS 05-34181, and by the Faculty of Arts and Sciences and the Department of Statistics of Rutgers University. We would like to thank the participants of the conference, the contributors to the volume, and the anonymous reviewers. Thanks are also due to DIMACS for providing conference facilities, and to the members of the staff and the many graduate students from the Department of Statistics for their tireless efforts to ensure the success of the conference. Last but not least, we would like to thank Ms. Pat Wolf for her patience and meticulous attention to all details in handling the papers in this volume.

Regina Liu, William Strawderman and Cun-Hui Zhang
December 15, 2006

# Dedication

This volume is dedicated to our dear colleague Yehuda Vardi, who passed away in January 2005.

Yehuda was born in 1946 in Haifa, Israel. He earned a B.S. in Mathematics from Hebrew University, Jerusalem, an M.S. in Operations Research from the Technion, Israel Institute of Technology and a Ph.D. under Jack Kiefer at Cornell University in 1977.

Yehuda served as a Scientist at AT&T's Bell Laboratories in Murray Hill before joining the Department of Statistics at Rutgers University in 1987. He served as the department chair from 1996 until he passed away. Yehuda was a dynamic and influential chair. He led the department with great energy and vision. He also provided much service to the statistical community by organizing many research conferences, and serving on the editorial boards of several statistical and engineering journals. He was an elected fellow of the *Institute of Mathematical Statistics* and *International Statistical Institute*. His research was supported by numerous grants from the *National Science Foundation* and other government agencies.

Yehuda was a leading statistician and a true champion for interdisciplinary research. He developed key algorithms which are now widely used for emission tomographic PET and SPECT scanners. In addition to his work on medical imaging, he coined the term "network tomography" in his pioneering paper on the problem of estimating source-destination traffic based on counts in individual links or "road sections" of a network. This problem has since blossomed into a full-fledged field of active research. His work on unbiased estimation based on biased data was a fundamental contribution in the field, and was recently rediscovered as a powerful general tool for the popular Markov chain Monte Carlo method. He has explored many other areas of statistics, including data depth and positive linear inverse problems with applications in signal recovery. His seminal contributions played a leading role in advancing the scientific fields in question, while enriching statistics with important applications.

Yehuda was not just a scientist with remarkable breadth and insight. He was also a wonderful colleague and friend, and a constant source of encouragement and humor. We miss him deeply.

Regina Liu, William Strawderman and Cun-Hui Zhang

# Contributors to this volume

Adhikari, A., *Avaya Labs*

Bickel, P. J., *University of California, Berkeley*

Cabrera, J., *Rutgers University*
Cao, J., *Bell Laboratories, Alcatel-Lucent Technologies*
Chang, I-S., *National Health Research Institutes*
Chen, A., *Bell Laboratories, Alcatel-Lucent Technologies*
Chien, L.-C., *National Health Research Institutes*
Cohen, A., *Rutgers University*

Denby, L., *Avaya Labs*

Eyheramendy, S., *Oxford University*

Fang, J., *Rutgers University*

Gijbels, I., *Katholieke Universiteit Leuven*

Hsiung, C. A., *National Health Research Institutes*

Jornsten, R., *Rutgers University*

Kolassa, J. E., *Rutgers University*

Landwehr, J. M., *Avaya Labs*
Lawrence, E., *Los Alamos National Laboratory*
Li, B., *Tsinghua University*
Lin, D., *University of North Carolina*
López-Pintado, S., *Universidad Pablo de Olavide*

Madigan, D., *Rutgers University*
Mallows, C., *Avaya Labs*
Mandel, M., *The Hebrew University of Jerusalem*
Meloche, J., *Avaya Labs*
Meyer, M., *University of Georgia*
Michailidis, G., *University of Michigan*

Nair, V. N., *University of Michigan*
Naus, J., *Rutgers University*

Pal, J. K., *University of Michigan*

Rinott, Y., *Hebrew University*

Sackrowitz, H., *Rutgers University*
Shepp, L., *Rutgers University*
Shlomo, N., *Southampton University*

Singh, K., *Rutgers University*
Strawderman, W. E., *Rutgers University*

Tang, W., *Rutgers University*

Vardi, Y., *Rutgers University*
Vitale, R. A., *University of Connecticut*

Wen, C.-C., *Tamkang University*
Woodroofe, M., *University of Michigan*
Wu, Y.-J., *Chung Yuan Christian University*

Xie, M., *Rutgers University*

Ying, Z., *Columbia University*
Yu, C.-R., *Rutgers University*

Zhang, C.-H., *Rutgers University*
Zhang, J., *Rutgers University*